\documentclass[12pt]{article}
\begin{document}
\begin{center}
\large{\bf{Discrepancy of $LS$-sequences of partitions}}
\bigskip

{\large{\bf{Ingrid Carbone}}}

Universit\`a della Calabria
\end{center}
\bigskip

\begin{abstract}
In this paper we give a precise estimate of the discrepancy of a class of uniformly distributed sequences of partitions. Among them we found a large class having low discrepancy (which means of order $ {\frac {1}{N}}$). One of them is the Kakutani-Fibonacci sequence.\end{abstract}

\section {\bf Introduction and preliminaries}

S. Kakutani studied in [K] the notion of {\it uniformly distributed sequences of partitions} of the interval $[0,1]$. He introduced the following construction.
Fix a number $\alpha \in \,]0,1[$. If $\pi$ is any partition of  $[0,1]$, its $\alpha$-refinement, denoted by $\alpha \pi$, is obtained subdividing the longest interval(s) of $\pi$ in proportion $\alpha/(1-\alpha)$. By $\alpha^n \pi$ we denote the $\alpha$-refinement of $\alpha^{n-1}\pi$.

Let $\omega=\{[0,1]\}$ be the trivial partition of $[0,1]$. The sequence $\{\alpha^n \omega\}$ will be called the Kakutani $\alpha$-sequence.
\bigskip

\noindent {\bf Definition 1.1} Given a sequence of partitions $\{\pi_n\}$ of $[0,1]$, with 

$$\pi_n = \{ [y_{i-1}^n , y_i^n[\,, 1 \le i \le k(n)\}\,,$$
  
  \smallskip
 \noindent we say that it is {\it uniformly distributed (u.d.)}  if for any continuous function $f$ on $[0,1]$ we have 

$$\lim_{n\rightarrow \infty} \frac{1}{k(n)}\sum_{i=1}^{k(n)} f(y_i^n)= \int_0^1 f(t)\,d t.\,\eqno(1)$$

\bigskip

We can now state Kakutani's result.
\bigskip

\noindent {\bf Theorem 1.2} {\it The sequence $\{\alpha^n \omega\}$ is uniformly distributed.}
\bigskip

This result caught the attention of several authors in the late seventies, when a different proof of the theorem was given ([AF]), and other papers were devoted to a stochastic version of it (see  [L1], [L2], [PvZ], [vZ]).

Recently, the procedure introduced by Kakutani has been generalized in several directions.

In [ChV] this notion has been extended to compact metric spaces.

In [V2] Kakutani's splitting procedure has been generalized, producing a new class of u.d. sequences of partitions as follows.

\bigskip

\noindent {\bf Definition 1.3} Consider any non trivial finite partition $\rho$ of $[0,1]$. The {\it $\rho$-refinement} of a partition
$\pi$ of $[0,1]$ (which will be denoted by $\rho \pi$) is obtained by
subdividing all the intervals of
$\pi$ having maximal length positively (or directly) homothetically to $\rho$.
\bigskip

Obviously, if $\rho =\{[0, \alpha[,[\alpha, 1]\}$, then the $\rho$-refinement is just
Kakutani's $\alpha$-refinement. 

As in Kakutani's case, we can iterate the splitting procedure. The $\rho$-refinement of $\rho \pi$ will be denoted by $\rho^2 \pi$ and, for all $n\in I\!\!N$, the  $\rho$-refinement of $\rho^{n-1} \pi$ will be denoted by $\rho^n \pi$.
If $\{\rho^n\omega\}$ denotes the sequence of successive $\rho$-refinements of the trivial partition $\omega$,  the following theorem holds (see [V2, Theorem 2.7]):
\bigskip

\noindent {\bf Theorem 1.4} {\it The sequence $\{\rho^n \omega\}$ is uniformly distributed.}

\bigskip

In [CV1] the splitting procedure has been generalized to higher dimensions, providing a sequence of nodes in the hypercube $[0,1]^d$ which is uniformly distributed. 

In [IV] the authors study uniform distribution on fractals.

In [CV2] the authors present a von Neumann type theorem (see [vN]), which provides uniformly distributed sequences of partitions of the interval $[0,1]$ starting from sequences of partitions whose diameter tends to 0 when $n\rightarrow \infty$.

In [DI] the authors give upper estimates of the discrepancy of $\rho$-refinements of the interval $[0,1]$.

Actually, the theory of uniformly distributed sequences of partition is deeply connected to the theory of uniformly distributed of sequences of points, which have a large application in quasi-Monte Carlo methods (see [N]).

For a large overview on uniformly distributed sequences of points, see [DT] and [KN].

Let us consider now a special class of $\rho$-refinements of the trivial partition $\omega$. 
\bigskip

\noindent {\bf Definition 1.5} Let us fix two positive integers $L$ and $S$ and let $0<\beta<1$ be the real number  such that  $ L\beta + S \beta^2=1$. Denote by $\rho\,_{L,S}$  the partition defined by $L$ ``long" intervals having length $\beta$ followed by S ``short" intervals having length $\beta^2$. 
By $\{\rho^n_{L,S}\,\omega\}$ (or $\{\rho^n_{L,S}\}$ for short) we denote the sequence of successive $ \rho_{L,S}$-refinements of the trivial partition $\omega$. They will be called $LS$-{\it sequences}.

\bigskip
 It is clear that  the partition $ \rho^n_{L,S}$ is obtained by dividing all the longest intervals of $ \rho^{n-1}_{L,S}$ homothetically with respect to $ \rho_{L,S}$, and that each partition $ \rho^n_{L,S}$  contains only two kinds of intervals: the long intervals have length $\beta^n$ while the short ones have length $\beta^{n+1}$. 

By Theorem 1.4 we know that  the sequence $\{\rho^n_{L,S}\}$ is uniformly distributed, but for this particular simple class of sequences  a direct proof was given in the lecture notes of a course on Uniform Distribution delivered by A. Vol\v{c}i\v{c} during the ``Workshop on Measure Theory and Real Analysis" held in 2002 ([V1]).
\bigskip

In this paper we give precise estimates of the discrepancy of $LS$-sequences.

Among the family of $LS$-sequences having low discrepancy, there is the Kakutani $\alpha$-refinement with $\alpha+\alpha^2=1$ (or, equivalently, $L=S=1$). Apart from the obvious case $\alpha=1/2$, this is the only Kalutani sequence for which the exact discrepancy is now known. This sequence may be called the Kakutani-Fibonacci sequence.

\section {\bf  $LS$-sequences of partitions}

This section will be devoted to the study of $LS$-sequences of partitions and their discrepancy.

Denote with $t_n$ the total number of intervals of $ \rho^n_{L,S}$, with $l_n$ the number of its long intervals and with $s_n$ the number of its short intervals. Obviously we have the following relations:

$$ t_n= s_n+l_n, \,\,\,\,\,\,\,\,\,\,l_n=L\,l_{n-1}+s_{n-1}, \,\,\,\,\,\,\,\,\,\,s_n=S\,l_{n-1}.$$

\bigskip

The terms of the sequence $t_n$ can be calculated solving the following difference equation
$$  t_n=L\, t_{n-1}+S\,t_{n-2}$$

\smallskip
\noindent with initial condition $t_0=1$ and $t_1=L+S$.

We easily get
$$ t_n= \frac{1+S\beta} {1+S \beta^2} \left(\frac {1} {\beta^n}\right)- \frac{S \beta-S \beta^2} {1+S \beta^2} (-S \beta)^n .\eqno (2)$$

\smallskip
We note that $l_n$ and $s_n$ satisfy the same difference equation, but the initial conditions are different since $l_0=1$, $l_1=L$ and $s_0=0$, $s_1=S$, respectively.
It is convenient to introduce the constants

$$A= \frac{1+S\beta} {1+S \beta^2}\,\,\,\,\,\,\,\,\rm{  and }\,\,\,\,\,\,\,\,B=\frac{S \beta-S \beta^2} {1+S\beta^2} ,$$

\bigskip

\noindent which will frequently appear in the sequel, so that formula (2) can be rewritten in the following way
$$t_n= \frac {A-B (-S \beta ^2)^n} {\beta^n}.\eqno (3)$$
\bigskip
Since $S\beta^2<1$, we see that $t_n$ has the same order as $1\over {\beta^n}$ when $n \rightarrow \infty.$

The following remark will play a role in identifying the different cases in the statement of Theorem 2.3.

\bigskip
\noindent {\bf Remark  2.1} From $L \beta+S \beta^2=1$ we get 

$$\beta= \frac{\sqrt{L^2+4S}-L}{2S},$$ 

\noindent from which we derive that

$$S\beta<1\,\, {\rm if\,\, and \,\,only\,\, if} \,\, S< L+1,$$
$$S\beta=1\,\, {\rm if\,\, and \,\,only\,\, if} \,\, S= L+1,$$
$$S\beta>1\,\, {\rm if\,\, and \,\,only\,\, if} \,\, S> L+1.$$

\bigskip
Now we recall the definition of {\it discrepancy} and {\it star-discrepancy}.
\bigskip

\noindent {\bf Definition 2.2}  Given a finite subset  $W=\{ w_1,w_2, \dots, w_N  \}$ of the interval $[0,1]$, the {\it discrepancy} of $W$ is defined as

$$ D(W) = \sup_{0 \le a < b \le 1}  \left | \frac {1} {N} \sum _{j=1}^{N} \chi _{[a,\, b[} (w_j) - (b-a)\right | ,$$

\bigskip
\noindent while the {\it star-discrepancy}  is defined as

$$ D^{\star}(W) = \sup_{0 < b \le 1}  \left | \frac {1} {N} \sum _{j=1}^{N} \chi _{[0, \,b[} (w_j) - b\,\right | .$$

\bigskip\bigskip

If we consider a sequence of partitions $\{\pi_n\}$ of the interval $[0, 1]$, with $\pi_n=  \{[y_{i-1}^n,  y_i^n[ , 1 \le i \le k(n)\}$, we may define 

$$ D(\pi_n) = \sup_{0 \le a < b \le 1}  \left | \frac {1} {k(n)} \sum _{j=1}^{k(n)} \chi _{[a, \,b[} (y_j^n) - (b-a)\right | $$

\noindent and
$$ D^{\star}(\pi_n) = \sup_{0 < b \le 1}  \left | \frac {1} {k(n)} \sum _{j=1}^{k(n)} \chi _{[0, \,b[} (y_j^n) - b\,\right | . \eqno (4)$$

\bigskip

It is well known that $\{ \pi_n  \}$ is uniformly distributed is and only if $D( \pi_n) \rightarrow 0$ when $n \rightarrow \infty$, and that $D^\star(W) \le D(W) \le 2 D^\star (W)$. 

The speed of the discrepancy is important in the applications.
It is also well known that the best discrepancy is $1 \over n$ and it is attained, for example, by  Knapowski's sequence $ \{ [ {{i-1 } \over {n}  }, {i \over n} [ , 1 \le i \le n\}$.

We are now ready to determine the star-discrepancy of the sequence  $ \{\rho_{L,S}^n \}$.

\bigskip

\noindent {\bf Theorem 2.3} {\it The discrepancy of the sequence of partitions  $ \{\rho_{L,S}^n \}$,  as  
$n \rightarrow \infty$, behaves in the following way:

\bigskip

i) $D( \rho_{L,S}^n  ) \sim {1 \over {t_n}}$   if  $S <L+1$,

\bigskip

ii) $D( \rho_{L,S}^n  ) \sim {{\log t_n} \over {t_n}}$  if  $S =L+1$,

\bigskip  

iii) $D( \rho_{L,S}^n ) \sim {1 \over {t_n}^ \gamma}$   if  $S> L+1$, where $ \gamma = 1+ \frac {\log(S \beta)} {\log \beta}$.}

\bigskip

\noindent {\it Proof.}  We denote the intervals of $ \rho^n_{L,S}$ as follows

$$ \rho^n_{L,S}=  \{[y_{i-1}^n,  y_i^n[ , 1 \le i \le t_n)\} $$

\smallskip
\noindent and we study the behavior of the star-discrepancy of $\rho_{L,S}^n$.

Denote by $L_p$ and $S_p$, respectively, a long interval and a short interval of the partition $\rho_{L,S}^p $, and first of all let us evaluate for $n\ge p$ the differences 

$$  \frac {1} {t_n} \sum _{j=1}^{t_n} \chi _{L_p}(y_j^n) - \lambda ( L_p) $$

\noindent and
$$  \frac {1} {t_n} \sum _{j=1}^{t_n} \chi _{S_p}(y_j^n) - \lambda ( S_p ) , $$

\smallskip
\noindent where of course $\lambda$  denotes the Lebesgue measure and $\lambda ( L_p ) = \beta^p$ and $\lambda ( S_p ) = \beta^{p+1}$.

We observe that, since $\rho_{L,S}^n$ is a refinement of $\rho_{L,S}^p$, each long interval $L_p$ of $\rho_{L,S}^p$ is the union of consecutive intervals of  $\rho_{L,S}^n $ (with $n \ge p$) and, since the splitting procedure is ``self-similar", this union reproduces up to a factor $\beta^{-p}$ the partition $\rho_{L,S}^{n-p}$. Therefore we have, for any $n\ge p$,

$$ \sum _{j=1}^{t_n} \chi _{L_p}(y_j^n) = t_{n-p} . $$
 
With a simple calculation we obtain, for $n\ge p$, 

$$  \frac {1} {t_n} \sum _{j=1}^{t_n} \chi _{L_p}(y_j^n) - \lambda ( L_s )  = \frac {t_{n-p}}{t_n} - \beta^p = $$
$$=-\frac {B}{t_n} \,\,(-S \beta)^n \left [(-S \beta)^{-p} - \beta ^p \right ] . \eqno (5)$$

\bigskip

If we consider now a short interval $S_p$ of $\rho_{L,S}^p  $, it becames a long one in  $\rho_{L,S}^{p+1}  $, which means that it is of the type $L_{p+1}$, so we can rewrite formula (5) in the following way:

$$  \frac {1} {t_n} \sum _{j=1}^{t_n} \chi _{S_p}(y_j^n) - \lambda ( S_p )  = \frac {t_{n-(p+1)}}{t_n} - \beta^{p+1} = $$
$$=-\frac {B}{t_n} \,\,(-S \beta)^n \left[(-S \beta)^{-p-1} - \beta ^{p+1} \right] \eqno (6)$$

\noindent for every $ n \ge p+1$.

Fix now $b\in\, ]0,1[$.

Let $[b_1^{n-1}, b_2^{n-1}[$ be the interval of $\rho_{L,S}^{n-1}$ containing $b$. Since the number of points of $\rho_{L,S}^{n-1}$ contained in $[b_1^{n-1}, b_2^{n-1}[$ is $1$ if this interval is a short one, and it is $L+S$ if it is a long one, we have

$$ \frac {1} {t_n} \sum _{j=1}^{t_n} \chi _{[0, \,b[} (y_j^n) - b \le \frac {1} {t_n} \sum _{j=1}^{t_n} \chi _{[0,\, b_2^{n-1}[} \,(y_j^n) - b_1^{n-1} = $$

$$=  \frac {1} {t_n} \sum _{j=1}^{t_n} \chi _{[0, \,b_1^{n-1}[} \,(y_j^n) - b_1^{n-1}+ \frac {1} {t_n} \sum _{j=1}^{t_n} \chi _{[b_1^{n-1}, \,b_2^{n-1}[} \,(y_j^n)\le$$

$$\le \frac {1} {t_n} \sum _{j=1}^{t_n} \chi _{[0, \,b_1^{n-1}[} \,(y_j^n) - b_1^{n-1}+ \frac {L+S} {t_n}\eqno (7)$$

\bigskip
\noindent and

$$ \frac {1} {t_n} \sum _{j=1}^{t_n} \chi _{[0, \,b[} (y_j^n) - b \ge \frac {1} {t_n} \sum _{j=1}^{t_n} \chi _{[0,\, b_1^{n-1}[}\, (y_j^n) - b_2^{n-1} = $$

$$=  \frac {1} {t_n} \sum _{j=1}^{t_n} \chi _{[0, \,b_2^{n-1}[} \,(y_j^n) - b_2^{n-1}- \frac {1} {t_n} \sum _{j=1}^{t_n} \chi _{[b_1^{n-1}, \,b_2^{n-1}[}\, (y_j^n)=$$

$$= \frac {1} {t_n} \sum _{j=1}^{t_n} \chi _{[0,\, b_2^{n-1}[}\, (y_j^n) - b_2^{n-1}- \frac {L+S} {t_n}.\eqno (8)$$

\bigskip
We will estimate (4) using the upper and the lower bounds in (7) and (8), so that we will evaluate 
$ \frac {1} {t_n} \sum _{j=1}^{t_n} \chi _{[0,\, b_k^{n-1}[}\, (y_j^n) - b_k^{n-1}$ for $k=1,2$. 
To this purpose it will be convenient to represent  $[0,\,b_k^{n-1}[$ (from now on we will write 
$[0,\,b^{n-1}[$  for short) as the union of consecutive intervals stemming from the partitions $\rho_{L,S}^i $ for $i \le n-1$.

Consider first all the consecutive intervals $I_1^1, I_2^1, \dots, I_{m_1}^1$ of $\rho_{L,S}^1$ such that  $ \bigcup_{i=1}^{m_1} I_i^1 \subset [0,\,b^{n-1}[$ ($m_1$ could be zero).

Next, take all the consecutive intervals $I_1^2, I_2^2, \dots, I_{m_2}^2$ of $\rho_{L,S}^2$ contained in  $ [0,\,b^{n-1}[\,\, \setminus  \bigcup_{i=1}^{m_1} I_i^1$ (again $m_2$ could be zero).

Proceed this way taking, at the final step,  all the $m_{n-1}$ (possibly zero) consecutive intervals $I_1^{n-1}, $ $I_2^{n-1}, \dots, I_{m_{n-1}}^{n-1}$ of $\rho_{L,S}^{n-1}  $ such that 
 
$$ \bigcup_{i=1}^{m_{n-1}} I_i^{n-1}  \,\,\,  \subset \,\,\, [0,\,b^{n-1}[\,\, \setminus   \, \bigcup_{p=1}^{n-2} \left ( \bigcup_{i=1}^{m_p} I_i^p \right),$$

\bigskip
\noindent so that at the end
 
$$  [0,\, b^{n-1}[ \,\,\, =\,\, \bigcup_{p=1}^{n-1}  \left( \bigcup_{i=1}^{m_p} I_i^p \right).       \eqno (9)$$ 

\bigskip
Thus $[0,\, b^{n-1}[\,$ is represented by the union of (at most) $n-1$ blocks of consecutive intervals of $\rho_{L,S}^p $ for $1 \le p \le n-1$.

Let $l^b_p$ and $s^b_p$ be, respectively,  the number of long intervals (denoted by  $L^p_i$) and short intervals (denoted by $S^p_i$), respectively, of $ \rho_{L,S}^p  $ contained in $[0,\,b[$ (obviously, $l_p^b+s_p^b=m_p$). Then using (5), (6) and (9), in order to evaluate the star-discrepancy we can write:

$$\frac {1} {t_n} \sum _{j=1}^{t_n} \chi _{[0,\,b^{n-1}[}  (y_j^n)  - b^{n-1}=  $$

$$=\frac {1} {t_n} \sum _{j=1}^{t_n} \chi _{ \left(\bigcup_{p=1}^{n-1}  \left( \bigcup_{i=1}^{m_p} I_i^p\right)\right)}  (y_j^n) - \lambda \left ( \bigcup_{p=1}^{n-1}  \left( \bigcup_{i=1}^{m_p} I_i^p \right ) \right) =$$
\smallskip
$$= \sum _{p=1}^{n-1}  \sum _{i=1}^{m_p}  \left(  \frac {1} {t_n} \sum _{j=1}^{t_n} \chi _{I_i^p}  (y_j^n)   - \lambda (I_i^p) \right)=$$
$$=  \sum _{p=1}^{n-1}  \sum _{i=1}^{l^b_p } \left(  \frac {1} {t_n} \sum _{j=1}^{t_n} \chi _{L_i^p}  (y_j^n)   - \lambda (L_i^p) \right)+ \sum _{p=1}^{n-1}  \sum _{i=1}^{s^b_p } \left(  \frac {1} {t_n} \sum _{j=1}^{t_n} \chi _{S_i^p}  (y_j^n)   - \lambda (S_i^p) \right)=$$
$$ = \frac {B}{\,t_n\,} (-S \beta)^n \sum _{p=1}^{n-1} \left \{ - l^b_p\, [(-S \beta)^{-p} - \beta ^p] -s^b_p \, [(-S \beta)^{-p-1} - \beta ^{p+1}] \right \} =$$
$$=  \frac {B}{\,t_n\,} (-S \beta)^n  \sum _{p=1,  p\,\, {\rm even}}^{n-1} \left[- l^b_p\, [(-S \beta)^{-p} - \beta ^p]  -s^b_p\, [(-S \beta)^{-p-1} - \beta ^{p+1} ]\right] +$$
$$ + \frac {B}{\,t_n\,} (-S \beta)^n  \sum _{p=1,  p\,\, {\rm odd}}^{n-1} \left[- l^b_p\,  [(-S \beta)^{-p} - \beta ^p]  -s^b_p\, [(-S \beta)^{-p-1} - \beta ^{p+1} ]\right]  = $$ 
$$=   \frac {B}{\,t_n\,} (-S \beta)^n \left\{  \sum _{p=1,  p\,\, {\rm even}}^{n-1} s^b_p\, [(S \beta)^{-p-1}  +\beta ^{p+1}]+ \sum _{p=1,  p\,\, {\rm odd}}^{n-1} l^b_p \,[(S \beta)^{-p}  +\beta ^{p}] \right\}-$$
$$- \frac {B}{\,t_n\,} (-S \beta)^n  \left\{ \sum _{p=1,  p\,\, {\rm even}}^{n-1} l^b_p\, [(S \beta)^{-p}  -\beta ^p]+ \sum _{p=1,  p\,\, {\rm odd}}^{n-1} s^b_p\, [(S \beta)^{-p-1}  -\beta ^{p+1}] \right\}.\eqno(10)$$

\bigskip
It is clear from the construction that $m_p \le L+S-1$, $0 \le l^b_p \le L$ and $0 \le s^b_p \le S$ for all $1 \le p \le n$, otherwise at least one of the intervals of $\bigcup_{i=1}^{m_p}I_i^p$ would be already present in $\bigcup_{i=1}^{m_{p-1}}I_i^{p-1}$. 

We have to distinguish the cases when $S\beta \not=1$ and $S\beta=1$. 

\bigskip
1) $S\beta \not=1$

We will give upper and lower estimates of (10) in the case $n$ even and in the case $n$ odd.

Assume first that $n$ is even and let $n=2h$, $h\ge1$. 

We note that the first term in (10) is positive, and the second term in (10) is negative since $B>0$ and $(S\beta)^{-p}+\beta^p=\frac {1-(S\beta^2)^p} {S\beta^2}$. Than with simple calculations we get the following upper and lower estimates:

$$-B (L+S-1) \frac {(S \beta)^{2h}}{t_{2h}}   \left\{   2 \left[  \frac {1-  (S\beta)^{2h} }  {(1- S^2\beta^2)(S \beta)^{2h-2}}  - \frac {1-\beta^{2h}}  {1-\beta^2}   \right] + \frac{1}{(S\beta)^{2h}} - \beta^{2h}    \right\} \le $$

$$\le \frac {1} {t_{2h}} \sum _{j=1}^{t_{2h}} \chi _{[0,\,b^{2h-1}[}  (y_j^{2h})  -  b^{2h-1} \le $$

$$\le B (L+S-1) \frac {(S \beta)^{2h}}{\,t_{2h}\,}  \left\{  2  \left[   \frac {1- (S \beta)^{2h-2}}  {(1- S^2\beta^2) (S \beta)^{2h-1} }   +
\beta   \frac {\beta^2-\beta^{2h}}  {1-\beta^2}      \right] + \frac {1} {S\beta} + \beta  \right\}.  \eqno (11)$$

\bigskip
Assume now that $n$ is odd and let $n=2h+1$, with $h\ge1$. Then with the same arguments from (10) we derive the following upper and lower bounds:

$$-B(L+S-1)\frac {(S \beta)^{2h+1}}{\,t_{2h+1}\,} \times$$
$$\times  \left\{   2 \left[  \frac {1-  (S\beta)^{2h-2} }  {(1- S^2\beta^2)(S \beta)^{2h-1}}  +\beta
  \frac {\beta^2-\beta^{2h}}  {1-\beta^2}  \right]     +\frac{1} {S\beta} + \beta +\frac{1}{(S \beta)^{2h+1}} + \beta^{2h+1}  \right\}  \le$$

$$\le \frac {1} {t_{2h+1}} \sum _{j=1}^{t_{2h+1}} \chi _{[0,\,b^{2h}[}  (y_j^{2h+1}) -  b^{2h} \le $$

 $$\le  2B (L+S-1)\frac { (S \beta)^{2h+1}}{\,t_{2h+1}\,}\left\{      \frac {1-  (S\beta)^{2h+2} }  {(1- S^2\beta^2)(S \beta)^{2h}}  - \frac {1-\beta^{2h+2}}  {1-\beta^2}    \right\}.  \eqno (12)$$

\bigskip
2) $S \beta=1$
 
 Equation (10) becomes

$$\frac {1} {t_n} \sum _{j=1}^{t_n} \chi _{[0,\,b^{n-1}[}  (y_j^n)  -  b^{n-1}  =$$
$$=  \frac {B}{\,t_n\,} (-1)^n \left\{  \sum _{p=1,  p\,\, {\rm even}}^{n-1} m_{p, S} [1  +\beta ^{p+1}]+ \sum _{p=1,  p\,\, {\rm odd}}^{n-1} m_{p, L} [1  +\beta ^{p}] \right\}-$$
$$- \frac {B}{\,t_n\,}  (-1)^n\left\{ \sum _{p=1,  p\,\, {\rm even}}^{n-1} m_{p, L} [1  -\beta ^p]+ \sum _{p=1,  p\,\, {\rm odd}}^{n-1} m_{p, S} [1 -\beta ^{p+1}] \right\} .\eqno (13)$$

\bigskip
If we assume $n=2h$, with $h\ge1$,  from equation (13) we derive the following upper and lower bounds

$$- \frac {B (L+S-1)}{\,t_{2h}\,} \left\{2h+1-2\,\frac{1 - \beta^{2h}}{1- \beta^2} -\beta^{2h}\right\} \le$$
$$ \le \frac {1} {t_{2h}} \sum _{j=1}^{t_{2h}} \chi _{[0,\,b^{2h-1}[}  (y_j^{2h}) - b^{2h-1} \le$$

$$\le   \frac {B(L+S-1) }{\,t_{2h}\,} \left\{  2h-1+2 \beta    \frac {\beta^2-\beta^{2h}}  {1-\beta^2}   +\beta \right\}. \eqno (14) $$ 

\bigskip 

If $n=2h+1$, with $h\ge1$, from (13) we derive the following inequalities

$$- \frac {B(L+S-1)}{\,t_{2h+1}\,}  \left\{  2h+ 2\beta    \frac {\beta^2-\beta^{2h}}  {1-\beta^2}  +\beta +\beta^{2h+1}\right\} \le$$

$$\le \frac {1} {t_{2h+1}} \sum _{j=1}^{t_{2h+1}} \chi _{[0,\,b^{2h}[}  (y_j^{2h+1}) - b^{2h}\le$$
  
$$ \le  \frac {B(L+S-1)}{\,t_{2h+1}\,} \left\{ 2h- 2\,\, \frac{\beta^{2} - \beta^{2h+2}}{1- \beta^2} \right\}. \eqno (15)$$ 

\bigskip 
Since we enclosed the generic $b \in [0,\,1[$ into the interval $[b_1^{n-1}, b_2^{n-1}[$ and we obtained the corresponding estimates (11), (12), (14) and (15), the same estimates hold for $b \in [0,\,1[$ and therefore they can be applied to evaluate the star discrepancy.

It is time now to distinguish among the three cases explicited in the statement of our theorem, in order to prove the various asymptotic behaviors of the star-discrepancy.
\bigskip

i) $S<L+1$

\bigskip
From Remark 2.1 this condition is equivalent to the case $S \beta <1$.

If $n$ is even, from (7), (8) and (11) simple calculations give 

$$-[L+S+ B  (L+S-1)] \times$$
$$\times\left\{ 2 \left[    \frac {S^2 \beta^2-  (S\beta)^{2h+2} }  {1- S^2\beta^2}  -(S \beta)^{2h} \frac {1-\beta^{{2h}}}  {1-\beta^2} \right]+1-(S\beta^2)^{2h}\right \} \le$$

$$\le  t_{2h} \left[\sup_{0 < b \le 1} \left( \frac {1} {t_{2h}} \sum _{j=1}^{t_{2h}} \chi _{[0,\,b[}  (y_j^{2h}) - b\right) \right] \le$$

\smallskip
$$\le L+S+B (L+S-1) \times $$
$$\times \left \{  2  \left[   \frac {S{\beta}- (S\beta)^{2h-1}}  {1- S^2\beta^2 }   +
\beta  (S \beta)^{2h} \frac {\beta^2-\beta^{2h}}  {1-\beta^2}      \right] + (S\beta)^{2h-1} +\beta (S \beta)^{2h} \right \} . \eqno (16)$$

\bigskip

If $n$ is odd, from (7), (8) and (12) 
we obtain
$$-[L+S+ B  (L+S-1)]  \left\{  2 \left[   \frac {S^2 \beta^2-  (S\beta)^{2h} }  {1- S^2\beta^2}   + \beta(S \beta)^{2h+1} \frac {\beta^2-\beta^{2h}}  {1-\beta^2} \right]\right\} -$$
$$ -[L+S+ B  (L+S-1)] \left\{(S\beta)^{2h} +\beta(S \beta)^{2h+1}+1+(S\beta^2)^{2h+1}\right\} \le$$
$$\le  t_{2h+1}\left[\sup_{0 < b \le 1} \left( \frac {1} {t_{2h+1}} \sum _{j=1}^{t_{2h+1}} \chi _{[0,\,b[}  (y_j^{2h+1}) - b\right) \right] \le$$

$$\le L+S+ B  (L+S-1) \left\{   2 \left[  \frac {S \beta-  (S\beta)^{2h+3} }  {1- S^2\beta^2}  -(S\beta)^{2h+1} \,  \frac {\,1-\beta^{2h+2}}  {1-\beta^2}  \right]      \right\}. \eqno (17)$$ 

\bigskip
\noindent Since $\beta<1$, $S\beta<1$ and $S\beta^2<1$, there exist the following limits:

$$  \lim_{h\rightarrow \infty} \left\{2  \left[   \frac {S{\beta}- (S\beta)^{2h-1}}  {1- S^2\beta^2 }   +
\beta  (S \beta)^{2h} \frac {\beta^2-\beta^{2h}}  {1-\beta^2}      \right] + (S\beta)^{2h-1} +\beta (S \beta)^{2h}\right\}= $$

$$= \lim_{h\rightarrow \infty}
  \left\{   2 \left[  \frac {S \beta-  (S\beta)^{2h+3} }  {1- S^2\beta^2}  -(S\beta)^{2h+1} 
  \frac {\,\,1-\beta^{2h+2}}  {1-\beta^2}  \right]      \right\}= \frac{2S \beta} {1-S^2\beta^2}$$

\noindent and

$$\lim_{h\rightarrow \infty} \{  2 \left[   \frac {S^2 \beta^2-  (S\beta)^{2h} }  {1- S^2\beta^2}   + \beta(S \beta)^{2h+1} \frac {\beta^2-\beta^{2h}}  {1-\beta^2} \right] +$$
\bigskip
$$+ (S\beta)^{2h} +\beta(S \beta)^{2h+1}+1+(S\beta^2)^{2h+1}\}=$$

$$=\lim_{h\rightarrow \infty} \left\{  2\left[    \frac {S^2 \beta^2-  (S\beta)^{2h+2} }  {1- S^2\beta^2}  -(S \beta)^{2h} \frac {1-\beta^{2h}}  {1-\beta^2}\right]+1-(S\beta^2)^{2h} \right\}=$$

$$=\frac {2S^2 \beta^2 }  {1- S^2\beta^2}+1,  
$$

\bigskip
\noindent so that the sequences of upper and lower bounds in (16) and (17) are bounded. Then  we conclude that there exists a positive constant $C_1$, independent on n, depending only on $L$ and $S$ (actually, $\beta$ depends on $L$ and $S$), such that  for all $n \in I\!\!N$

 $$  t_n\sup_{0 < b \le 1} \left|  \frac {1} {t_n} \sum _{j=1}^{t_n} \chi _{[0,\,b[}  (y_j^n) - b\right| \le C_1,$$
 
 \bigskip
 \noindent which implies that, as $n \rightarrow \infty$,
 
 $$D^*( \{\rho_{L,S}^n \} ) \sim {1 \over {t_n}}.$$

\bigskip
 
 ii) $S =L+1$
 
 \bigskip
First of all we note that, by Remark 2.1, we have $S \beta =1$. 

If $n$ is even, from (7), (8) and (14)  we obtain 

$$- \frac {B (L+S-1) }{\log t_{2h}} \left[{2h}+1-2\frac{1 - \beta^2}{1- \beta^{2h}} -\beta^{2h}\right]\le$$ 

$$\le \frac{t_{2h}} {\log t_{2h}}\left[   \sup_{0 < b \le 1} \left( \frac {1} {t_{2h}} \sum _{j=1}^{t_{2h}} \chi _{[0,\,b[}  (y_j^{2h}) - b\right) \right] \le $$
\bigskip
$$\le   \frac {B(L+S-1) }{\log t_{2h}\,} \left\{  2h-1+2 \beta    \frac {\beta^2-\beta^{2h}}  {1-\beta^2}    +\beta \right\}. \eqno (18) $$

\bigskip
If $n$ is odd, we obtain from (7), (8) and (15) that
\bigskip

$$-\frac {B (L+S-1)}{\log t_{2h+1}} \left\{  2h+ 2\beta    \frac {\beta^2-\beta^{2h}}  {1-\beta^2}  +\beta +\beta^{2h+1}\right\}\le$$

$$\le \frac{t_{2h+1}} {\log t_{2h+1}}\left[\sup_{0 < b \le 1} \left(  \frac {1} {t_{2h+1}} \sum _{j=1}^{t_{2h+1}} \chi _{[0,\,b[}  (y_j^{2h+1}) - b \right)\right] \le $$
\bigskip
$$\le \frac {B (L+S-1)}{\log t_{2h+1} } \left\{ 2h- 2\,\, \frac{\beta^{2} - \beta^{2h+2}}{1- \beta^2} \right\}. \eqno (19)$$

\bigskip
From formula (3) we have

$$ \log t_n = \log \left(\frac{A}{\beta^n}-B(-S\beta^n)\right)=$$

$$= \log \left( \frac{A-B(-S\beta^2)^n} {\beta^n}\right) = \log(A-B(-S\beta^2)^n) - n \log \beta$$

\bigskip
\noindent for all  $n \in I\!\!N$ so that, since $\beta<1$ and $S\beta^2<1$,

$$ \lim_{h\rightarrow \infty} \frac {2h-1+2 \beta    \frac {\beta^2-\beta^{2h}}  {1-\beta^2}    +\beta}{\log t_{2h}}=\lim_{h\rightarrow \infty} \frac {2h+1-2\frac{1 - \beta^2}{1- \beta^{2h}} -\beta^{2h}}{\log t_{2h}}=$$

$$ =\lim_{h\rightarrow \infty} \frac {2h- 2\,\, \frac{\beta^{2} - \beta^{2h+2}}{1- \beta^2} }{\log t_{2h+1}}= \lim_{h\rightarrow \infty}  \frac { 2h+ 2\beta    \frac {\beta^2-\beta^{2h}}  {1-\beta^2}  +\beta +\beta^{2h}}{\log t_{2h+1}} =\frac{ -1}{\log\beta}.  $$

\bigskip

Taking the above limits into account, we conclude from (18) and (19) that there exists a positive constant $C_2$, independent on n, depending only on $L$ and $S$, such that  
 
 $$ \frac{t_n}{\log t_n} \sup_{0 < b \le 1} \left|\frac {1} {t_n} \sum _{j=1}^{t_n} \chi _{[0,\,b[}  (y_j^n) - b \right| \le C_2,$$
 
 \bigskip
 \noindent which implies that, as  $n \rightarrow \infty$, 
 
$$D^*( \{\rho_{L,S}^n \} ) \sim {{\log t_n} \over {t_n}}.$$

\bigskip
iii) $S>L+1$

\bigskip

By Remark 2.1 this condition is equivalent to $S \beta >1$.

If $n$ is even, from  (7), (8) and (11) we have:

 $$ -  \frac { B(L+S-1)}{\beta^{2h}t_{2h}}   \left\{  2 \left[  \frac {1-  (S\beta)^{{2h}} }  {(1- S^2\beta^2)(S \beta)^{{2h}-2}}  - \frac {1-\beta^{{2h}}}  {1-\beta^2}   \right] + \frac{1}{(S\beta)^{2h}} - \beta^{2h}      \right\} \le$$

$$\le  \frac{1} {(S\beta^2)^{2h}} \left[\sup_{0 < b \le 1}  \left( \frac {1} {t_{2h}} \sum _{j=1}^{t_{2h}} \chi _{[0,\,b[}  (y_j^{2h}) - b \right)\right]\le  $$

$$\le   \frac {B(L+S-1)}{\beta^{2h}t_{2h}} \left\{  2  \left[   \frac {1- (S \beta)^{{2h}-2}}  {(1- S^2\beta^2) (S \beta)^{{2h}-1} }   +
\beta   \frac {\beta^2-\beta^{2h}}  {1-\beta^2}      \right] + \frac {1} {S\beta} + \beta  \right\}. \eqno (20) $$
 
\bigskip
If $n$ is odd, from (7), (8) and (12)  we get

$$-\frac {B(L+S-1)}{\beta^{2h+1}t_{2h+1}} \times$$

$$\times \left\{   2 \left[  \frac {1-  (S\beta)^{2h-2} }  {(1- S^2\beta^2)(S \beta)^{2h-1}}  +\beta
  \frac {\beta^2-\beta^{2h}}  {1-\beta^2}  \right]   +\frac{1} {S\beta} + \beta +\frac{1}{(S \beta)^{2h+1}} + \beta^{2h+1} \right\} \le$$

$$\le \frac{1} {(S\beta^2)^{2h+1}}\left[\sup_{0 < b \le 1}\left(\frac {1} {t_{2h+1}} \sum _{j=1}^{t_{2h+1}} \chi _{[0,\,b[}  (y_j^{2h+1}) - b\right)\right] \le$$

\bigskip
$$\le \frac {2B(L+S-1)}{\beta^{2h+1}t_{2h+1}}  \left[  \frac {1-  (S\beta)^{2h+2 } } {(1- S^2\beta^2)(S \beta)^{2h} } -
  \frac {1-\beta^{2h+2}}  {1-\beta^2}  \right] .\eqno (21)   $$

 \bigskip
 Since $S\beta>1$, $S\beta^2<1$ and $\beta<1$, using formula (3) we have

$$ \lim_{h\rightarrow \infty} \frac{1} {\beta^{2h}t_{2h}} \left\{  2  \left[   \frac {1- (S \beta)^{2h-2}}  {(1- S^2\beta^2) (S \beta)^{2h-1} }   +
\beta   \frac {\beta^2-\beta^{2h}}  {1-\beta^2}      \right] + \frac {1} {S\beta} + \beta  \right\}= $$ 

$$= \lim_{h\rightarrow \infty} \frac{1} {A-B(-S\beta^2)^{2h}}\times$$

$$\times \left\{  2  \left[   \frac {1- (S \beta)^{2h-2}}  {(1- S^2\beta^2) (S \beta)^{2h-1} }   +
\beta   \frac {\beta^2-\beta^{2h}}  {1-\beta^2}      \right] + \frac {1} {S\beta} +\beta  \right\}= $$ 

$$= \lim_{h\rightarrow \infty}  \frac {1}{\beta^{2h+1}t_{2h+1}}   \Big\{  2 \left[  \frac {1-  (S\beta)^{2h-2} }  {(1- S^2\beta^2)(S \beta)^{2h-1}}  +\beta
  \frac {\beta^2-\beta^{2h-1}}  {1-\beta^2}  \right]     +$$

$$+\frac{1} {S\beta} + \beta +\frac{1}{(S \beta)^{2h+1}} + \beta^{2h+1} \}=$$

$$= \frac {1}{A}  \left\{  2 \left[\frac{1}{S\beta(S^2\beta^2-1)}+
  \frac {\beta^3}  {1-\beta^2 }\right]   +\frac{1} {S\beta} + \beta \right\} $$

 \bigskip
\noindent and 
 
$$  \lim_{h\rightarrow \infty}  \frac {1}{\beta^{2h}t_{2h}}  \left\{  2 \left[  \frac {1-  (S\beta)^{{2h}} }  {(1- S^2\beta^2)(S \beta)^{{2h}-2}}  - \frac {1-\beta^{{2h}}}  {1-\beta^2}   \right] + \frac{1}{(S\beta)^{2h}} - \beta^{2h}      \right\}=$$  
  
 $$= \lim_{h\rightarrow \infty}  \frac {2}{\beta^{2h+1}t_{2h+1}}    \left[  \frac {1-  (S\beta)^{2h+2 }}  {(1- S^2\beta^2)(S \beta)^{2h}}  -
 \frac {1-\beta^{2h+2}}  {1-\beta^2}  \right]      = $$

 $$= \frac {2}{A} \left\{  \frac{S^2\beta^2}{S^2\beta^2-1}-\frac {1}  {1-\beta^2}    \right\} . $$

\bigskip
From (20), (21) and the previous limits, we conclude that there exists a positive constant $C_3$, independent on n, depending only on $L$ and $S$, such that  

 $$ \frac{1}{(S\beta^2)^n} \sup_{0 < b \le 1} \left| \frac {1} {t_n} \sum _{j=1}^{t_n} \chi _{[0,\,b[}  (y_j^n) - b\right| \le C_3$$
 
 \bigskip
 \noindent for all $n \in I\!\!N$, which implies that,  as  $n \rightarrow \infty$,

 $$D( \{\rho_{L,S}^n \} ) \sim {(S\beta^2)^n}$$
 
 \bigskip
 
 It remains only to observe how  $(S\beta^2)^n$ can be written in terms of $t_n$. 
 
In fact,  $(S\beta^2)^n=\beta^\gamma$, where $ \gamma = 1+ \frac {\log(S \beta)} {\log \beta}<1$ since 
 $S\beta^2=(S\beta)\beta=\beta^{1+c}$  with $ c=\frac{\log(S\beta)}{\log\beta}<0$.
Consequently, $(S\beta^2)^n$ and ${1 \over {t_n}^ \gamma}$ have the same order at infinity since

 $$ \frac{(S\beta^2)^n} {1/(t_n)^\gamma}= \frac{(S\beta^2)^n} {(\beta^n)^\gamma}\left( A+B(-S\beta^2)^n\right)^\gamma=\frac{(\beta^\gamma)^n} {(\beta^n)^\gamma}\left( A+B(-S\beta^2)^n\right)^\gamma\rightarrow A^\gamma$$
 \noindent as  $n \rightarrow \infty$, and finally, as $n \rightarrow \infty$,  we have
 
 $$D^*( \{\rho_{L,S}^n \} ) \sim {1 \over {t_n}^ \gamma},$$

 \noindent where $ \gamma = 1+ \frac {\log(S \beta)} {\log \beta}$.
 
 Since the discrepancy and the star-discrepancy are equivalent, the theorem is proved. $\star$

\section {\bf Conclusions} 

Uniformly distributed sequences of partitions having low discrepancy turn out to be useful in quasi-Monte Carlo methods, as it is clear from formula (1). 

$LS$-sequences of partitions offer the advantage that each partition is the refinement of the previous one.
Nevertheless, the disadvantage of these sequences  lies in the fact that the function has to be evaluated only on "blocks" of points, so that it is not possible to use intermediate number of points between  $t_n$ and $t_{n+1}$.

There exists a connection between u.d. sequences of partitions and u.d. sequences of points. In [V2, Theorem 3.4] it has been proved that a random reordering of the points of a u.d. sequence of partitions gives with probability one a u.d. sequence of points. That theorem, even if interesting from a theoretical point of view, does not give any information about the discrepancy of these random sequences of points.

The author conjectures that it is possible to associate to $LS$-sequences having low discrepancy sequences of points having low discrepancy too, i.e. of order $\log N\over N$. Some encouraging results have been already obtained, and will be presented in a forthcoming paper.

\bigskip\bigskip
\noindent {\bf Acknowledgements}

\bigskip
The author wishes to  express her gratitude to Aljo\v{s}a Vol\v{c}i\v{c} for his useful suggestions and remarks, and for his critical reading of the manuscript. 

\bigskip \bigskip

\noindent {\bf References}
 \bigskip
 
\noindent [A-F] R. L. Adler, L. Flatto, Uniform distribution of Kakutani's intervals splitting procedure, {\it Z. Wahrschleinlichkeitstheorie verw. Gebiete {\bf38}} (1977) 253-259. 
\smallskip\smallskip

 \noindent  [Ch-V] F. Chersi,  A. Vol\v{c}i\v{c},  $\lambda$-equidistributed sequences of partitions and a theorem of the de Bruijn-Post type,  {\it Annali Mat. Pura Appl. (4) {\bf 162 } } (1992) 23-32.
\smallskip\smallskip

\noindent  [C-V1] I. Carbone, A. Vol\v{c}i\v{c}, Kakutani splitting procedure in
higher dimension, {\it  Rend. Ist. Matem. Univ. Trieste} {\bf 39}  (2007) 119-126.
\smallskip\smallskip

\noindent  [C-V2] I. Carbone, A. Vol\v{c}i\v{c}, A von Neumann theorem for uniformly distributed sequences of partitions, {\it submitted}.
\smallskip\smallskip

\noindent [D-I] M. Drmota, M. Infusino,   On the discrepancy of some generalized Kakutani's sequences of partitions,  {\it preprint}. 
\smallskip\smallskip

\noindent [D-T] M. Drmota,  R. F. Tichy,  {\it Sequences, discrepancies and applications},  Lecture Notes in Mathematics {\bf 1651}, Springer Verlag, Berlin, 1997. 
\smallskip\smallskip

\noindent [I-V] M. Infusino, A. Vol\v{c}i\v{c}, Uniform distribution on fractals, {\it Uniform Distribution Theory} {\bf 4} (2009), n. 2, 47-58.  
\smallskip\smallskip

\noindent [K] S. Kakutani, A problem on equidistribution on the unit interval $[0,1]$, {\it  Measure theory (Proc. Conf., Oberwolfach, 1975)},  pp. 369--375. {\it Lecture Notes in Math.} {\bf 541}, Springer, Berlin, 1976.
\smallskip\smallskip

\noindent [K-N] L. Kuipers, H. Niderreiter, {\it Uniform distribution of sequences.  Pure and Applied Matematics}. Wiley-Interscience, New York - London-Sidney, 1974.
\smallskip\smallskip

\noindent [L1] J. C. Lootgieter, Sur la r\'epartition des suites de Kakutani, {\it C. R. Acad. Sci. Paris S\'er. AB}  {\bf \,285} (1977), no. 5, A403-A406.
\smallskip\smallskip

\noindent [L2] J. C. Lootgieter, Sur la r\'epartition des suites de Kakutani, {\it C. R. Acad. Sci. Paris SŽr. AB} {\bf 286} (1978), no. 10, A459-A461
\smallskip\smallskip

\noindent [N] H. Niderreiter, {\it Random number generation and quasi-Monte Carlo Methods}, CBMS-NSF Regional Conference Series in Applied Math., 1992.
\smallskip\smallskip

\noindent [P-vZ]  R. Pyke, W. R. van Zwet, Weak convergence results for the Kakutani interval splitting procedure, {\it Ann. Probability} {\bf \,32} (2004), no. 1A, 380-423.
\smallskip\smallskip

\noindent [vZ] W.R.  van Zwet, A proof of Kakutani's conjecture on random subdivision of longest intervals, {\it Ann. Probability} {\bf 6} (1978), no. 1, 133-137.
\smallskip\smallskip

\noindent [V1] A. Vol\v{c}i\v{c},  Uniformly distributed sequences of partitions. {\it Workshop on Measure Theory and Real Analysis, Mondello (PA), Italy}, July 8-17 (2003).
\smallskip\smallskip

\noindent [V2] A. Vol\v{c}i\v{c},  A generalization of Kakutani's splitting procedure, {\it to appear in Annali Mat. Pura e Appl.}

\smallskip\smallskip

\noindent [vN]  J. von Neumann, Gleichm\"assig dichte Zahlenfolgen, {\it Mat.
Fiz. Lapok} {\bf 32}  (1925) 32-40.
\smallskip\smallskip

\noindent [W] H. Weyl, \"Uber ein Problem aus dem Gebiete der diophantischen Approximationen, {\it Nach. Ges. Wiss. G\"ottingen, Math.-phys. Kl.} (1914), 234-244.

\end{document}